\documentclass{eptcs}

\usepackage{mathtools}
\usepackage{amssymb}
\usepackage{nicefrac}
\usepackage{enumitem}
\usepackage{fancybox}
\usepackage{environ}
\usepackage{wrapfig}
\usepackage{pdfrender}
\usepackage{xspace}
\usepackage{pict2e}
\usepackage{shadow}
\usepackage{xcolor}

\newtheorem{defin}{Definition}[section] % [chapter]
\newtheorem{lemma}[defin]{Lemma}

\newtheorem{definition}[defin]{Definition}
\newtheorem{example}[defin]{Example}
\newtheorem{notation}[defin]{Notation}
\newtheorem{convention}[defin]{Convention}

\ifdefined\mathsl\else
  \DeclareMathAlphabet{\mathsl}{\encodingdefault}{\rmdefault}{\mddefault}{\sldefault}
  \SetMathAlphabet{\mathsl}{bold}{\encodingdefault}{\rmdefault}{\bfdefault}{\sldefault}
\fi

\def\QEDbox{\def\sboxsep{.3em}\def\sdim{.1em}\raisebox{.3em}{\shabox{}}}
\def\QED{\hfill\QEDbox}

\newcommand{\after}{\mathrel{\circ}}
\newcommand{\idmap}[1][]{\ensuremath{\mathsl{id}_{#1}}}
\newcommand{\tuple}[1]{\langle#1\rangle}
\newcommand{\ket}[1]{\ensuremath{|{\kern.1em}#1{\kern.1em}\rangle}}
\newcommand{\bigket}[1]{\ensuremath{\big|{\kern.1em}#1{\kern.1em}\big\rangle}}
\newcommand{\ketstrut}{\vrule height 8.5pt depth 4.5pt width 0pt}
\newcommand{\Bigket}[1]{\ensuremath{\left|\ketstrut{\kern.1em}\right.{\kern-.2em}#1{\kern-.2em}\left.\ketstrut{\kern0em}\right>}}
\newcommand{\shortplus}{\ensuremath{{\kern-1.5pt}+{\kern-1.5pt}}}
\newcommand{\Mlt}{\ensuremath{\mathcal{M}}}
\newcommand{\natMlt}{\ensuremath{\mathcal{N}}}
\newcommand{\Lst}{\ensuremath{\mathcal{L}}}

\newcommand{\zero}{\ensuremath{\mathbf{0}}}
\newcommand{\acc}{\ensuremath{\mathsl{acc}}}
\newcommand{\head}{\ensuremath{\mathsl{hd}}}
\newcommand{\delete}{\ensuremath{\mathsl{del}}}
\newcommand{\quota}{\ensuremath{\mathsl{qu}}}
\newcommand{\concat}{\ensuremath{\mathbin{+{\kern-.5ex}+}}}
\newcommand{\unit}{\ensuremath{\mathsl{unit}}}
\newcommand{\flatten}{\ensuremath{\mathsl{flat}}}

\begin{document}

\title{Counting Votes with Multisets\thanks{Version of \today}}

%\author{Bart Jacobs, Michael Johnson, Richard Buckland}

\author{Bart Jacobs
%\institute{NICTA\\ Sydney, Australia}
\institute{Radboud University 
\\
%University of New South Wales\thanks{A fine university.}\\
Nijmegen, The Netherlands}
\email{bart@cs.ru.nl}
\and
Michael Johnson \qquad\qquad Richard Buckland
\institute{University of New South Wales \\
Sydney, Australia}
\email{michael.johnson1@unsw.edu.au \qquad richardb@unsw.edu.au}
}

%% \institute{iHub, 
%%  Radboud  University Nijmegen, The Netherlands.
%% }  
%% \email{bart@cs.ru.nl}

%\date{\small\today}

\def\titlerunning{Counting Votes with Multisets}
\def\authorrunning{Jacobs, Johnson, Buckland}

\maketitle

\begin{abstract}
A multiset is a `set' in which elements may occur multiple
times. These structures are ideal for expressing the outcome of an
election, for instance of the form 60 `yes' and 40 `no'. Moreover,
multisets are a useful datatype in vote counting algorithms. This will
be illustrated in three different forms of vote counting, known as:
`instant-runoff', `De Borda', and `single transferrable vote'. The
relevant abstract properties of multisets are: (1)~they form a (free)
commutative monoid, and (2)~they form a functor, and (3)~also a
monad. This paper illustrates how such categorical properties can be
put to good use in deriving and expressing election outcomes. The
emphasis is not on the (elementary) category theory involved, but on
its application in voting systems.
\end{abstract}

% Keywords: vote counting, multiset, functor, monad

\section{Introduction}

Counting votes seems obvious. Suppose there are two candidates $A$ and
$B$, and one hundred people have voted, precisely once, either $A$ or
$B$. Then we can arrive at an outcome of, say, 72 votes for $A$ and 28
votes for $B$. This paper looks at the datatypes that are appropriate
for counting votes and for reporting the vote outcome. For instance,
at the end of the voting period one may have a \emph{list} of length
100, containing 72 $A$'s and 28 $B$'s, in some order. The vote
outcome, we claim, is best described as a \emph{multiset}. In this
case it will be written as $72\ket{A} + 28\ket{B}$. In such a multiset
the order of the elements is irrelevant, only their multiplicity
counts. These multiplicities are written as numbers before `ket'
symbols $\ket{\cdot}$, with the candidates written inside the ket.

What we have described is a very simple election, with only two
candidates, and only one vote per candidate. There may be more
candidates or multiple options, and voters may express preferences, in
the form of a list of candidates, say in descending order of
preference. The term `preferendum' is sometimes used for such
elections. There may be good reasons to organise a preferendum, with
multiple weighted options, instead of a referendum, with only two
options.  A referendum may be hijacked to express discontent on other
matters --- as perhaps happened in the case of Brexit. A preferendum,
in contrast, invites a more nuanced stance, and seems less susceptible
to external influence. These political considerations are interesting,
but are not the topic of this paper, and only serve as motivation.

It turns out that it is not trivial to describe, explain, or implement
these vote counting mechanisms with multiple weighted options. Things
become even more difficult when surplus votes (after reaching a
threshold) are transferred, in some form. The two main points of this
paper are the following.
\begin{enumerate}
\item Multisets are the appropriate datatype for vote counting and
  reporting.

\item Some elementary techniques from category theory capture the
  operations on multisets that are relevant in vote couting algorithms
  (like instant-runoff, De Borda, single transferable votes). In
  particular, the fact that multisets form a functor, a monad, and
  carry a (free) monoid structure turns out to be remarkably useful.
\end{enumerate}

\noindent The paper aims to demonstrate these points, mainly via
examples of different types of elections, with different forms of vote
counting. The (partly implicit) point that the paper aims to make is
that describing these vote counting algorithms in terms of multisets
makes them `abstract' so that they can be described in a few lines.
One can debate whether abstractness contributes to proper
understanding, but mathematicians generally think so. And many
computer scientists think so too, since it makes it easier for them to
produce provably correct implementations. Again largely implicit, is
the motivation that vote counting should happen in a transparant and
correct manner.

The category theory used in this paper is relatively elementary, but
it is extremely useful. The contribution of this paper is thus not
categorical, but lies in recognising and exploiting the relevant
categorical structure in these applications, especially for multisets.

This paper explains different vote count methods mainly via examples,
with relatively small numbers. The calculations can, in principle, be
checked by hand. For this paper we have used Python scripts, from (an
updated version of) the EfProb library~\cite{ChoJ17}.

The paper devotes quite a bit of space to explaining multisets and
their categorical structure. These explanations are aimed at
non-category theorists, to make the material accessible outside the
circle of (categorical) experts. Readers who wish to learn more
category theory are referred to the extensive
literature~\cite{Awodey06,BarrW90,Leinster14,Pierce91,MacLane71,Simons11,Cheng22,Perrone24b}.
The next two sections contain a gentle introduction to the datatypes
of multisets and lists. Lists are ubiquitous in computer science and
mathematics, but multisets seem to be an orphaned, under-valued
datatype. Section~\ref{InstantSec} describes the `functoriality' of
multisets and puts this property to good use in describing so-called
instant-runoff voting to determine a winner from lists of
preferences. The fact that multisets form a monad is used in
Section~\ref{BordaSec} for De Borda vote counting, with different
weights. Finally, Section~\ref{STVSec} describes the mechanism of
single transferable votes for multi-member electorates. As will be
discussed in those sections, these different types of elections are
actually used, at various places in the world. Notably Australia has a
rich tradition of (subtle) vote counting. The paper ends with a number
of conclusions.

\section{A first look at multisets}\label{MltIntroSec}

A set is a collection of elements. For instance, one may have a set of
colours $S = \{R, B, G\}$, for red, blue and green. In a set, elements
occur at most once. There is the set of natural numbers $\mathbb{N} =
\{0,1,2,3,\ldots\}$ in which each number occurs precisely once. The
elements in a set are not ordered. We can also write $\mathbb{N} =
\{1,0,3,2,5,4,\ldots\}$.

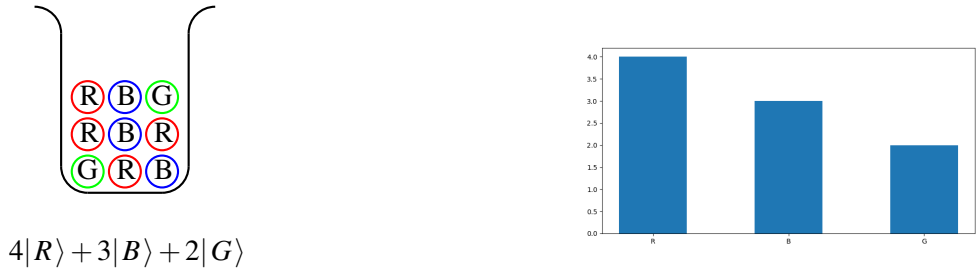
\begin{figure}
\[ \begin{array}{c}
\vcenter{\hbox{\begin{picture}(70,70)
\thicklines
% urn
\put(68, 60){\oval(20, 20)[tl]}
\put(0, 60){\oval(20, 20)[tr]}
\put(20, 10){\oval(20, 20)[bl]}
\put(48, 10){\oval(20, 20)[br]}
\put(20, 0){\line(1, 0){28}}
\put(10, 10){\line(0, 1){50}}
\put(58, 10){\line(0, 1){50}}
% with balls
\color{green}
\put(20,8){\circle{12}}
\color{black}
\put(16,5){G}
\color{red}
\put(34,8){\circle{12}}
\color{black}
\put(31,5){R}
\color{blue}
\put(48,8){\circle{12}}
\color{black}
\put(45,5){B}
\color{red}
\put(20,22){\circle{12}}
\color{black}
\put(17,19){R}
\color{blue}
\put(34,22){\circle{12}}
\color{black}
\put(31,19){B}
\color{red}
\put(48,22){\circle{12}}
\color{black}
\put(45,19){R}
\color{red}
\put(20,36){\circle{12}}
\color{black}
\put(17,33){R}
\color{blue}
\put(34,36){\circle{12}}
\color{black}
\put(31,33){B}
\color{green}
\put(48,36){\circle{12}}
\color{black}
\put(44,33){G}
\end{picture}}}
\\
\\
4\ket{R} + 3\ket{B} + 2\ket{G}
\end{array}
\hspace*{10em}
\vcenter{\hbox{\includegraphics[scale=0.25]{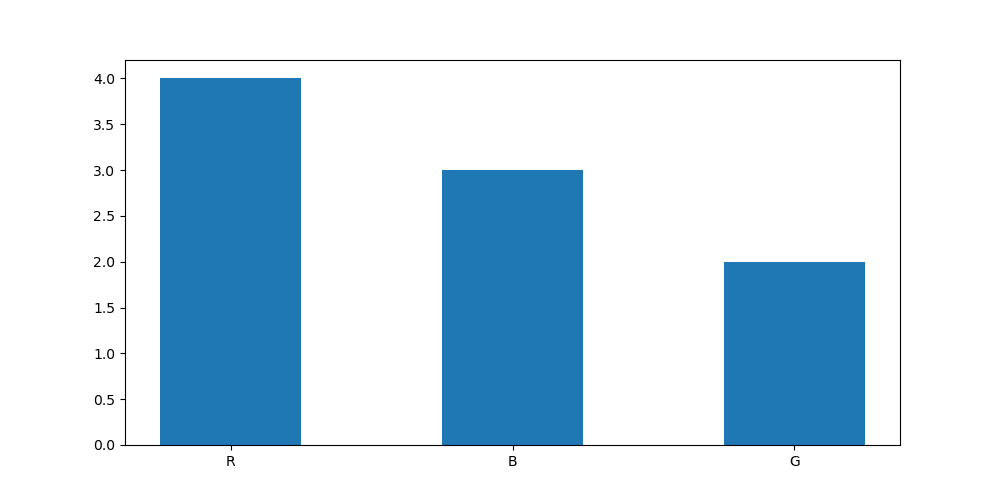}}}
\]
\caption{An urn with four red, three blue and two green balls,
  described as multiset on the left, with the corresponding bar chart
  on the right.}
\label{UrnChartFig}
\end{figure}

There are situations where one wishes to be able to talk about
multiple occurrences of elements. Figure~\ref{UrnChartFig} shows an
urn / vase containing coloured balls. One sees that the urn contains
nine balls in total, of which four are red, three blue and two green.
This is a typical example of a multiset. We shall use special `ket'
notation and write this urn / bag / multiset as:
\[ 4\ket{R} + 3\ket{B} + 2\ket{G} 
\qquad\mbox{ expressing }\qquad
\left\{\begin{array}{l}
\mbox{4 red}
\\
\mbox{3 blue}
\\
\mbox{2 green.}
\end{array}
\right. \]

\noindent The verticle bar $\mid$ and the right angle $\rangle$ together
form what it is called a ket $\ket{\cdot}$. Inside such a ket we write
an element, like a color $R$ / $B$ / $G$. Before the ket we write how
often that element occurs in a multiset. These kets have no mathematical
meaning but are used as notation to separate the elements (inside)
from their multiplicities (upfront).

The important thing to note is that the urn in
Figure~\ref{UrnChartFig} may also be seen as a \emph{ballot box},
containing nine votes, namely four for option~$R$, three for
option~$B$, and two for option~$G$. It is at this stage unclear ---
and irrelevant --- what these three options are. The important point
is that options that voters express typically occur multiple times,
and thus can be described by multisets. The bar chart on the
right in Figure~\ref{UrnChartFig} tabulates the votes and clearly
shows that option~$R$ has most votes.

\begin{convention}
\label{MltConv}
Several conventions apply when we write multisets via kets. 
\begin{enumerate}
\item \label{MltConvOrder} The order in which the items are listed
  does not matter.  Thus we have equalities of multisets:
\[ \begin{array}{rcccl}
4\ket{R} + 3\ket{B} + 2\ket{G} 
& = &
2\ket{G} + 3\ket{B} + 4\ket{R}
& = &
3\ket{B} + 2\ket{G} + 4\ket{R}.
\end{array} \]

\item \label{MltConvZero} Zero occurrences are typically ommitted. If
  we write $Y$ for yellow, then:
\[ \begin{array}{rcl}
0\ket{Y} + 4\ket{R} + 3\ket{B} + 2\ket{G}
& = &
4\ket{R} + 3\ket{B} + 2\ket{G}.
\end{array} \]

\noindent We may write the term $0\ket{Y}$ to emphasise that there are
no (zero) yellow balls (or votes) in the urn, but usually this is not
so relevant.

\item \label{MltConvSum} The multiplicities of multiply occurring kets
  with the same element are added, as in:
\begin{equation}
\label{MltAddEqn}
\begin{array}{rcl}
2\ket{R} + 2\ket{G} + 3\ket{B} + 2\ket{R}
& = &
4\ket{R} + 3\ket{B} + 2\ket{G}.
\end{array}
\end{equation}

\noindent This makes sense, both in terms of urns and ballot boxes.
\end{enumerate}

\noindent This third equation will turn out to be very useful.
\end{convention}

\begin{figure}
\[ \vcenter{\hbox{\begin{picture}(270,90)
\thicklines
% urn1
\put(68, 60){\oval(20, 20)[tl]}
\put(0, 60){\oval(20, 20)[tr]}
\put(20, 10){\oval(20, 20)[bl]}
\put(48, 10){\oval(20, 20)[br]}
\put(20, 0){\line(1, 0){28}}
\put(10, 10){\line(0, 1){50}}
\put(58, 10){\line(0, 1){50}}
% with balls
\color{green}
\put(20,8){\circle{12}}
\color{black}
\put(16,5){G}
\color{red}
\put(34,8){\circle{12}}
\color{black}
\put(31,5){R}
\color{blue}
\put(48,8){\circle{12}}
\color{black}
\put(45,5){B}
\color{red}
\put(20,22){\circle{12}}
\color{black}
\put(17,19){R}
\color{blue}
\put(34,22){\circle{12}}
\color{black}
\put(31,19){B}
\color{red}
\put(48,22){\circle{12}}
\color{black}
\put(45,19){R}
\color{red}
\put(20,36){\circle{12}}
\color{black}
\put(17,33){R}
\color{blue}
\put(34,36){\circle{12}}
\color{black}
\put(31,33){B}
\color{green}
\put(48,36){\circle{12}}
\color{black}
\put(44,33){G}
% urn2
\put(168, 60){\oval(20, 20)[tl]}
\put(100, 60){\oval(20, 20)[tr]}
\put(120, 10){\oval(20, 20)[bl]}
\put(148, 10){\oval(20, 20)[br]}
\put(120, 0){\line(1, 0){28}}
\put(110, 10){\line(0, 1){50}}
\put(158, 10){\line(0, 1){50}}
% urn3
\put(268, 60){\oval(20, 20)[tl]}
\put(200, 60){\oval(20, 20)[tr]}
\put(220, 10){\oval(20, 20)[bl]}
\put(248, 10){\oval(20, 20)[br]}
\put(220, 0){\line(1, 0){28}}
\put(210, 10){\line(0, 1){50}}
\put(258, 10){\line(0, 1){50}}
% with balls
\color{blue}
\put(220,8){\circle{12}}
\color{black}
\put(216,5){B}
\color{red}
\put(234,8){\circle{12}}
\color{black}
\put(231,5){R}
\color{blue}
\put(248,8){\circle{12}}
\color{black}
\put(245,5){B}
\color{green}
\put(220,22){\circle{12}}
\color{black}
\put(217,19){G}
\color{green}
\put(234,22){\circle{12}}
\color{black}
\put(231,19){G}
\color{green}
\put(248,22){\circle{12}}
\color{black}
\put(245,19){G}
% bent arrows
\cbezier (35,50) (40,120) (122,120) (127,50)
\put(127,50){\vector(1,-5){2}}
\cbezier (235,50) (230,120) (145,120) (140,50)
\put(140,50){\vector(-1,-5){2}}
\end{picture}}}
\]
\caption{Addition of multisets can be represented via a join of
  contents of two urns / ballot boxes. The resulting urn in the middle
  will contain the sum of the two outer multisets, namely:
   $5\ket{R} + 5\ket{B} + 5\ket{G}$.}
\label{UrnAddFig}
\end{figure}
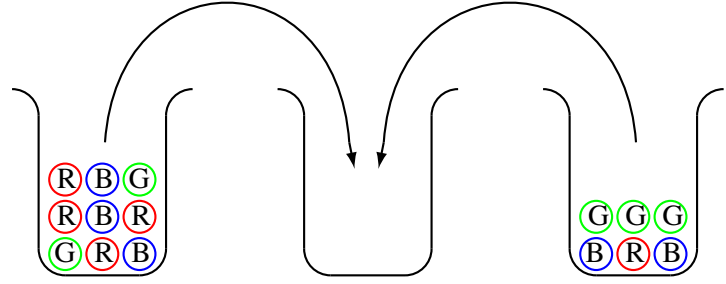

We are used to adding numbers, like in $2+3=5$. One can also add
multisets. This can be represented by pooring the contents of two urns
/ ballot boxes into a new urn, see Figure~\ref{UrnAddFig}. We write
this addition of multisets also as $+$. This may be confusing at
first, but is actually convenient. For instance the adding of
multisets in Figure~\ref{UrnAddFig} can be written as a mathematical
equation:
\[ \begin{array}{rcl}
\Big(4\ket{R} + 3\ket{B} + 2\ket{G} \Big) + 
   \Big(1\ket{R} + 2\ket{B} + 3\ket{G} \Big)
& = &
5\ket{R} + 5\ket{B} + 5\ket{G}.
\end{array} \]

\noindent In fact, this equation can be seen as a consequence of
Convention~\ref{MltConv}~\eqref{MltConvSum}, in the style of
Equation~\eqref{MltAddEqn}.

\begin{notation}
\label{NatMltNot}
For a set $X$, we write $\natMlt(X)$ for the set of multisets, written
as $n_{1}\ket{x_1} + \cdots + n_{K}\ket{x_K}$, with elements $x_{i}\in X$
and with natural numbers as multiplicities $n_{i}\in\mathbb{N}$. 

For a such a multiset $\varphi = \sum_{i} n_{i}\ket{x_i} \in \natMlt(X)$
we write $\|\varphi\| \in \mathbb{N}$ for the size of the multiset,
given by the sum of multiplicities:
\[ \begin{array}{rcccl}
\|\varphi\|
& = &
\Big\|\,\sum_{i} n_{i}\ket{x_i} \,\Big\|
& = &
\sum_{i} n_{i}.
\end{array} \]

For two multisets $\varphi,\psi\in\natMlt(X)$ we write $\varphi + \psi
\in\natMlt(X)$ for the sum of multisets, obtained by addition of
multiplicties, basically as in~\eqref{MltAddEqn}. This addition is
commutative, satisfying $\varphi + \psi = \psi + \varphi$ for all
multisets $\varphi, \psi\in \natMlt(X)$.

Writing $\zero \in\natMlt(X)$ for the empty / zero multisets, with
multiplicity zero for each element from $X$, we get $\varphi = \zero +
\varphi = \varphi + \zero$.
\end{notation}

In Figure~\ref{UrnAddFig}, the two urns on the left and right have
size~$9$ and~$6$. Their sum has size~$15$. Indeed, in general one has
$\big\|\varphi+\psi\big\| = \big\|\varphi\big\| + \big\|\psi\big\|$.
The general, mathematical statement is that $\natMlt(X)$ is a (free)
commutative monoid, and that the size map $\|-\| \colon \natMlt(X)
\rightarrow \mathbb{N}$ is a homomorphism of monoids, that preserves
addition (including zero). For (much) more information,
see~\cite{Jacobs26a}.

Now that we know about addition of multisets, we can already compute
outcomes of elections. We illustrate this in a simple preferendum.

\begin{example}
\label{ThreeOptionEx}
Consider an election with a set of four options $V = \{a,b,c,d\}$.
Voters can express a vote of the form $(x,y,z)$, where $x\in V$ is
their first choice that counts three times, $y\in V$ is their second
choice that counts twice, and $z$ is their third choice that counts
only once. In a valid vote $(x,y,z)$ these choices $x,y,z$ should all
be different. There are in this simple example seven votes, of the
form:
\[ (d,a,c) \quad (a,b,d) \quad (c,a,b) \quad (d,a,c) \quad (c,b,a) \quad 
   (a,b,c) \quad (d,b,c). \]

\noindent We can now calculate the vote outcome as a multiset over
$V$, that is, as an element of $\natMlt(V)$. We proceed by turning
each vote $(x,y,z)$ into a multiset $3\ket{x} + 2\ket{y} + 1\ket{z}$.
Then we can add all these multisets. This yields:
\[ \begin{array}{rcl}
\lefteqn{\Big(3\ket{d} + 2\ket{a} + 1\ket{c}\Big) +
   \Big(3\ket{a} + 2\ket{b} + 1\ket{d}\Big) +
   \Big(3\ket{c} + 2\ket{a} + 1\ket{b}\Big) +
   \Big(3\ket{d} + 2\ket{a} + 1\ket{c}\Big)}
\\[+0.3em]
& & \quad + \;
   \Big(3\ket{c} + 2\ket{b} + 1\ket{a}\Big) +
   \Big(3\ket{a} + 2\ket{b} + 1\ket{c}\Big) +
   \Big(3\ket{d} + 2\ket{b} + 1\ket{c}\Big)
\\[+0.3em]
& = &
13\ket{a} + 9\ket{b} + 10\ket{c} + 10\ket{d}.
\end{array} \hspace*{6em} \]

\noindent Thus, option~$a$ is most popular. This is an instance of a
De Borda count~\cite{Borda81,Emerson13,Emerson16}, which we shall
describe more systematically in Section~\ref{BordaSec}. The
calculations require some care when done by hand, but when you use a
programming language (or package) that supports multisets, the
counting can be done systematically, at a high level of abstraction.
\end{example}

\section{From lists to multisets}\label{ListSec}

One can imagine a situation where people have to write their vote, say
$Y$ for yes and $N$ for no, on a piece of paper and then hand it over
to some trusted person who counts the votes. This person ends up with
a list / pile of papers. The outcome can then be extracted by turning
the list into a multiset. We call this process \emph{accumulation} and
use a function $\acc$ to do this. This could look as follows.
\begin{equation}
\label{AccExEqn}
\begin{array}{rcl}
\acc\Big(\tuple{Y, Y, N, Y, N, N, N, Y, Y, Y}\Big)
& = &
6\ket{Y} + 4\ket{N}.
\end{array}
\end{equation}

\noindent Using lists instead of multisets as datatype for votes is
not such a good idea, mainly because it may leak information and thus
compromise confidentiality. Indeed, if someone has registered the
order in which people have voted, then seeing the election outcome as
a list of votes, makes it possible to reconstruct individual votes.

An additional advantage of using multisets is that it is easy to see
the outcome, as on the right in~\eqref{AccExEqn}. Seeing the list does
not immediately tell us the conclusion.

There is some mathematical structure at work that is nice to make
explicit.

\begin{notation}
\label{LstNot}
For a set $X$ we write $\Lst(X)$ for the set of (finite) lists with
elements from the set $X$. A typical element of $\Lst(X)$ is of the
form $\tuple{x_{1}, \ldots, x_{n}}$, with elements $x_{i} \in X$.  The
order of elements does matter in a list. Elements may occur multiple
times in a list.

Given two lists $\ell_{1}, \ell_{2} \in \Lst(X)$ we write $\ell_{1} \concat
\ell_{2} \in \Lst(X)$ for their concatenation. It is simply obtained
by gluiing them together, in order, as in:
\[ \begin{array}{rcl}
\tuple{x_{1}, \ldots, x_{n}} \concat \tuple{y_{1}, \ldots, y_{m}}
& = &
\tuple{x_{1}, \ldots, x_{n}, y_{1}, \ldots, y_{m}}.
\end{array} \]

\noindent Together with the empty list $\tuple{} \in \Lst(X)$ this
turns the set $\Lst(X)$ into a monoid --- but not a commutative one.

We write $\acc \colon \Lst(X) \rightarrow \natMlt(X)$ for the accumulation
function that turns a list into a multiset by counting occurrences
and forgetting the order. This may be written as:
\[ \begin{array}{rcl}
\acc\Big(\tuple{x_{1}, \ldots, x_{n}}\Big)
& = &
1\ket{x_1} + \cdots + 1\ket{x_n}.
\end{array} \]

\noindent This is a convenient and clever definition, since the
multiplicities of multiply occurring elements are automatically added
by Convention~\ref{MltConv}~\eqref{MltConvSum}. Mathematically, this
function $\acc$ does the vote counting, by turning a list of ballots
with expressed (single) votes into a multiset with counted votes.
\end{notation}

A basic fact is that accumulation maps concatenations to sums, as in:
\begin{equation}
\label{AccHomEqn}
\begin{array}{rclcrcl}
\acc\big(\ell_{1} \concat \ell_{2}\big)
& = &
\acc\big(\ell_{1}\big) + \acc\big(\ell_{2}\big)
& \qquad\mbox{ and }\qquad
\acc\big(\tuple{}\big)
& = &
\zero.
\end{array}
\end{equation}

\noindent Formulated more abstractly, the function $\acc \colon
\Lst(X) \rightarrow \natMlt(X)$ is a homomorphism of monoids.

This is directly relevant for voting. Consider a polling place with
three ballot boxes $\mathsl{BB}_{1}$, $\mathsl{BB}_{1}$,
$\mathsl{BB}_{3}$ in which lists $\ell_{1}$, $\ell_{2}$, $\ell_{3}$ of
votes are collected, and counted at the end of the day:
\[
\raisebox{-1.8em}{\rotatebox{90}{\ovalbox{\text{ voters }}}}
\hspace*{0.2em}
\left.\begin{array}{r}
\stackrel{\ell_1}{\longrightarrow} \mathsl{BB}_{1}
\\
\stackrel{\ell_2}{\longrightarrow} \mathsl{BB}_{2}
\\
\stackrel{\ell_3}{\longrightarrow} \mathsl{BB}_{3}
\end{array}\right\}
\hspace*{0.5em}\xrightarrow{\text{ vote counting }}\hspace*{0.5em}
\acc\big(\ell_{1}\big) + \acc\big(\ell_{2}\big) + \acc\big(\ell_{3}\big).
\]

\noindent The sum of the three accumulated lists on the
right-hand-side is the outcome of the election. The sum $+$ is a sum
of multisets. Notice that by using this sum $+$ of multisets, we
assume that each ballot box is counted separately. One could also
first put the piles of votes together via concatenation $\concat$ and
then count them all via accumulation. This gives as election outcome
$\acc\big(\ell_{1} \concat \ell_{2} \concat \ell_{3}\big)$.
Via~\eqref{AccHomEqn} we can be sure that the outcome is the same.

\section{Functoriality of multisets and instant-runoff voting}\label{InstantSec}

Instant-runoff voting is a vote count mechanism that is
used\footnote{according to Wikipedia
\href{https://en.wikipedia.org/wiki/Instant-runoff_voting}{en.wikipedia.org/wiki/Instant-runoff\_voting},
consulted Nov.~29, 2025} to elect members of the Australian House of
Representatives and the National Parliament of Papua New Guinea. It is
also used to elect the head of state in India, Ireland, and Sri
Lanka. The vote counting happens in several rounds, in which one
candidate with the least highest-preference-votes, is eleminated, but
where the other votes (with this eliminated candidate first) are still
included, see below for details.

This vote counting process is non-trivial, both in complexity and in
amount of work. It turns out that an abstract description can be given
in terms of multisets, using that multisets form a
\emph{functor}. This means that the multiset operation $\natMlt$ does
not only apply to sets, like in Notation~\ref{NatMltNot}, but also to
functions between them. This property is also called
\emph{functoriality}.

\begin{definition}
\label{LstMltFunDef}
Let $X,Y$ be two arbitrary sets, with a function $f\colon X
\rightarrow Y$ between them. In this situation we can form the two sets
$\Lst(X)$ and $\Lst(Y)$ of lists over $X$ and over $Y$, and also the
sets $\natMlt(X)$ and $\natMlt(Y)$ of multisets over $X$ and over $Y$.
\begin{enumerate}
\item One can define a function between these two sets of lists,
written as $\Lst(f) \colon \Lst(X) \rightarrow \Lst(Y)$, namely:
\begin{equation}
\label{LstFunEqn}
\begin{array}{rcl}
\textstyle\Lst(f)\Big(\tuple{x_{1}, \ldots, x_{n}}\Big)
& = &
\tuple{f(x_{1}), \ldots, f(x_{n})}.
\end{array}
\end{equation}

\item We can also define a function between the two sets of multisets,
  written as $\natMlt(f) \colon \natMlt(X) \rightarrow \natMlt(Y)$,
  via:
\begin{equation}
\label{MltFunEqn}
\begin{array}{rcl}
\textstyle\natMlt(f)\Big(\sum_{i} n_{i}\bigket{x_i}\Big)
& = &
\sum_{i} n_{i}\bigket{f(x_{i})}.
\end{array}
\end{equation}
\end{enumerate}
\end{definition}

The function $\Lst(f)$ in~\eqref{LstFunEqn} performs what functional
programmers call map-list. It simply applies the function $f$ to all
the elements of a list. The function $\natMlt(f)$ in~\eqref{MltFunEqn}
works similarly, but the sum property of
Convention~\ref{MltConv}~\eqref{MltConvSum} may kick in, when
different elements in $X$ are mapped by the function $f$ to the same
element in $Y$. We shall exploit this property below.

In terms of urns filled with coloured balls, the function $\natMlt(f)$
can be seen as a repainting of the balls in the urn, where the
function $f\colon X \rightarrow Y$ maps one set of colours to another
set.  For instance, let's write $Y$ for yellow and $P$ for purple, and
consider a repainting function $f \colon \{R,B,G\} \rightarrow
\{Y,P\}$, given by $f(R) = f(G) = Y$ and $f(B) = P$. Then we can
repaint the (contents) of the urn in Figure~\ref{UrnChartFig} via the
function $\natMlt(f)$, as elaborated in:
\[ \begin{array}{rcl}
\natMlt(f)\Big(4\bigket{R} + 3\bigket{B} + 2\bigket{G}\Big)
& \smash{\stackrel{\eqref{MltFunEqn}}{=}} &
4\bigket{f(R)} + 3\bigket{f(B)} + 2\bigket{f(G)}
\\[+0.2em]
& = &
4\bigket{Y} + 3\bigket{P} + 2\bigket{Y}
\\[+0.2em]
& = &
6\bigket{Y} + 3\bigket{P}.
\end{array} \]

\noindent There is no such thing as repainting of ballots in a ballot
box, but nevertheless we shall soon see the usefulness of
functoriality in vote counting. But first we collect (without proof)
some basic properties.

\begin{lemma}
\label{LstMltFunLem}
In the context of Definition~\ref{LstMltFunDef},
\begin{enumerate}
\item $\Lst$ preserves function composition and identities: $\Lst(g
  \after f) = \Lst(g) \after \Lst(f)$ and $\Lst(\idmap) = \idmap$;

\item $\natMlt$ also preserves composition and identities: $\natMlt(g
  \after f) = \natMlt(g) \after \natMlt(f)$ and $\natMlt(\idmap) = \idmap$;

\item $\Lst(f)$ preserves length and $\natMlt(f)$ preserves size; the
  latter can be expressed as $\big\|\natMlt(f)(\varphi)\big\| =
  \big\|\varphi\big\|$ for each multiset $\varphi\in\natMlt(X)$;

\item Accumulation is natural: $\acc \after \Lst(f) = \natMlt(f)
  \after \acc$, for each function $f$. \QED
\end{enumerate}
\end{lemma}

We shall describe instant-runoff vote counting via an example. It
involves two auxilary functions on lists, namely a `head' and `delete'
function $\head$ and $\delete_{a}$ for an element $a$. The head
function $\head$ selects the first element of a (non-empty) list, and
the delete function $\delete_{a}$ removes the element $a$ from the
entire list; it returns the list that remains.  For instance:
\[ \begin{array}{rclcrcl}
\head\Big(\tuple{c,d,a,b}\Big)
& = &
c
& \qquad\qquad &
\delete_{a}\Big(\tuple{c,d,a,b}\Big)
& = &
\tuple{c,d,b}.
\end{array} \]

We assume we have an election where one candidate must be chosen out
of set $X = \{a,b,c,d\}$. There are 100 voters and they each give a
list of preferences, $\tuple{x_{1}, x_{2}, x_{3}, x_{4}}$ for
$x_{i}\in X$, all different. We assume the 100 votes are divided in the
following way as top-down preference sequences.
\begin{equation}
\label{PrefSeqTable}
\hbox{\begin{tabular}{c||c|c|c|c|c|c|c|c|c|c|c|}
 & 7 & 2 & 7 & 8 & 7 & 19 & 12 & 12 & 8 & 6 & 12
\\
\hline
\hline
$1^{\mathsl{st}}$ & 
   $a$ & $a$ & $a$ & $b$ & $b$ & $b$ & $c$ & $d$ & $d$ & $d$ & $d$
\\
$2^{\mathsl{nd}}$ & 
   $b$ & $c$ & $d$ & $a$ & $c$ & $d$ & $a$ & $a$ & $a$ & $b$ & $c$
\\
$3^{\mathsl{rd}}$ & 
   $d$ & $d$ & $c$ & $c$ & $d$ & $c$ & $d$ & $b$ & $c$ & $c$ & $a$
\\
$4^{\mathsl{th}}$ & 
   $c$ & $b$ & $b$ & $d$ & $a$ & $a$ & $b$ & $c$ & $b$ & $a$ & $b$
\end{tabular}}
\end{equation}

\noindent We reformulate this table as a multiset
$\varphi_{0}\in\natMlt\big(\Lst(X)\big)$ over lists, of the form:
\[ \begin{array}{rcl}
\varphi_{0}
& = &
7\bigket{a, b, d, c} + 
   2\bigket{a, c, d, b} + 
   7\bigket{a, d, c, b} + 
   8\bigket{b, a, c, d}
\\[+0.2em]
& & \quad + \,
   7\bigket{b, c, d, a} + 
   19\bigket{b, d, c, a} + 
   12\bigket{c, a, d, b} + 
   12\bigket{d, a, b, c}
\\[+0.2em]
& & \quad + \,
   8\bigket{d, a, c, b} + 
   6\bigket{d, b, c, a} + 
   12\bigket{d, c, a, b}.
\end{array} \]

%% 7|['a', 'b', 'd', 'c']> + 
%%    2|['a', 'c', 'd', 'b']> + 
%%    7|['a', 'd', 'c', 'b']> + 
%%    8|['b', 'a', 'c', 'd']> + 
%%    7|['b', 'c', 'd', 'a']> + 
%%    19|['b', 'd', 'c', 'a']> + 
%%    12|['c', 'a', 'd', 'b']> + 
%%    12|['d', 'a', 'b', 'c']> + 
%%    8|['d', 'a', 'c', 'b']> + 
%%    6|['d', 'b', 'c', 'a']> + 
%%    12|['d', 'c', 'a', 'b']>

In an instant-runoff counting process one first looks at the top row
and one counts which candidates occurr there, with which
multiplicities. By inspecting the above table~\eqref{PrefSeqTable} we
see:
\begin{equation}
\label{ByHandEqn}
\left\{\begin{array}{lrcl}
\mbox{for $a$, } & 16 & = & 7 + 2 + 7 
\\
\mbox{for $b$, } & 34 & = & 8 + 7 + 19 
\\
\mbox{for $c$, } & 12 & & 
\\
\mbox{for $d$, } & 38 & = & 12 + 8 + 6 + 12.
\end{array} \right.
\end{equation}

\noindent Interestingly, these numbers can also be obtained via
functoriality of $\natMlt$, using the head function:
\[ \begin{array}{rcl}
\natMlt\big(\head\big)(\varphi_{0})
& \smash{\stackrel{\eqref{MltFunEqn}}{=}} &
7\bigket{\head(a, b, d, c)} + 
   2\bigket{\head(a, c, d, b)} + 
   7\bigket{\head(a, d, c, b)} + 
   8\bigket{\head(b, a, c, d)}
\\[+0.2em]
& & \quad + \,
   7\bigket{\head(b, c, d, a)} + 
   19\bigket{\head(b, d, c, a)} + 
   12\bigket{\head(c, a, d, b)} + 
   12\bigket{\head(d, a, b, c)}
\\[+0.2em]
& & \quad + \,
   8\bigket{\head(d, a, c, b)} + 
   6\bigket{\head(d, b, c, a)} + 
   12\bigket{\head(d, c, a, b)}
\\[+0.2em]
& = &
7\bigket{a} + 
   2\bigket{a} + 
   7\bigket{a} + 
   8\bigket{b} +
   7\bigket{b} + 
   19\bigket{b} + 
   12\bigket{c} + 
   12\bigket{d} +
   8\bigket{d} + 
   6\bigket{d} + 
   12\bigket{d}
\\[+0.2em]
& = &
16\bigket{a} + 
   34\bigket{b} + 
   12\bigket{c} + 
   38\bigket{d}
\end{array} \]

\noindent These numbers correspond to the above count by
hand~\eqref{ByHandEqn}. They are obtained via
Convention~\ref{MltConv}~\eqref{MltConvSum}.

We see two things: (1)~no-one has an absolute majority, so there is no
winner yet; (2)~candidate $c$ has the lowest number of first
preferences. With the instant-runoff method, this candidate $c$ is
removed as possible winner.

This removal again makes use of functoriality of $\natMlt$, but now
with $\delete_{c}$ as function. We elaborate the details of the
computation. It gives us a new multiset $\varphi_{1}$ with which
the process is repeated.
\[ \begin{array}{rcl}
\lefteqn{\varphi_{1} 
\hspace*{\arraycolsep}\coloneqq\hspace*{\arraycolsep}
\natMlt\big(\delete_{c}\big)(\varphi_{0})}
\\[+0.3em]
& \smash{\stackrel{\eqref{MltFunEqn}}{=}} &
7\bigket{\delete_{c}(a, b, d, c)} + 
   2\bigket{\delete_{c}(a, c, d, b)} + 
   7\bigket{\delete_{c}(a, d, c, b)} + 
   8\bigket{\delete_{c}(b, a, c, d)}
\\[+0.2em]
& & \quad + \,
   7\bigket{\delete_{c}(b, c, d, a)} + 
   19\bigket{\delete_{c}(b, d, c, a)} + 
   12\bigket{\delete_{c}(c, a, d, b)} + 
   12\bigket{\delete_{c}(d, a, b, c)}
\\[+0.2em]
& & \quad + \,
   8\bigket{\delete_{c}(d, a, c, b)} + 
   6\bigket{\delete_{c}(d, b, c, a)} + 
   12\bigket{\delete_{c}(d, c, a, b)}
\\[+0.2em]
& = &
7\bigket{a, b, d} + 
   2\bigket{a, d, b} + 
   7\bigket{a, d, b} + 
   8\bigket{b, a, d} +
   7\bigket{b, d, a} + 
   19\bigket{b, d, a}
\\[+0.2em]
& & \quad + \,
   12\bigket{a, d, b} + 
   12\bigket{d, a, b} +
   8\bigket{d, a, b} + 
   6\bigket{d, b, a} + 
   12\bigket{d, a, b}.
\\[+0.2em]
& = &
7\bigket{a, b, d} + 
   21\bigket{a, d, b} + 
   8\bigket{b, a, d} + 
   26\bigket{b, d, a} + 
   32\bigket{d, a, b} + 
   6\bigket{d, b, a}.
\end{array} \]

%% 7|['a', 'b', 'd']> + 
%%    21|['a', 'd', 'b']> + 
%%    8|['b', 'a', 'd']> + 
%%    26|['b', 'd', 'a']> + 
%%    32|['d', 'a', 'b']> + 
%%    6|['d', 'b', 'a']>

We can now iterate and start the second round, in essentially the same
way. At first positions we now get:
\[ \begin{array}{rcl}
\natMlt\big(\head\big)(\varphi_{1})
& \smash{\stackrel{\eqref{MltFunEqn}}{=}} &
28\bigket{a} + 34\bigket{b} + 38\bigket{d}.
\end{array} \]

Still no-one has an absolute majority, and $a$ has the lowest
number. Hence this candidate $a$ is deleted in this second step,
giving as next multiset:
\[ \begin{array}{rcccl}
\varphi_{2}
& \coloneqq &
\natMlt\big(\delete_{a}\big)(\varphi_{1})
& \smash{\stackrel{\eqref{MltFunEqn}}{=}} &
41\bigket{b, d} + 59\bigket{d, b}.
\end{array} \]

\noindent At this stage things begin to become clearer. We still start the third
round by computing the first preferences:
\[ \begin{array}{rcl}
\natMlt\big(\head\big)(\varphi_{2})
& \smash{\stackrel{\eqref{MltFunEqn}}{=}} &
41\bigket{b} + 59\bigket{d}.
\end{array} \]

\noindent We now have an absolute majority for $d$ and can declare
that this candidate is the winner. We see how functoriality of $\natMlt$
does the work for us. 

%% \bigskip

%% \textbf{To be added}: 
%% \begin{itemize}
%% \item Allow lists of preferences of arbitrary length. Then we will
%%   have to define head as a partial function $\head \colon \Lst(X)
%%   \rightarrow X \cup \{\bot\}$ and carry the undefined element $\bot$
%%   along in the calculations. This does is not a fundamental change.

%% \item Re-evaluate actual voting results in this manner, using the
%%   Python scripts that were also used in the above examples.

%% \item Write the actual vote count in pseudo code as an iteration
%%   involving the $\natMlt(\head)$ and $\natMlt(\delete_{x})$
%%   calculations.

%% \item Are there homomorphism properties for this instant-runoff, like
%%   in~\eqref{AccHomEqn}?
%% \end{itemize}

\section{De Borda counts}\label{BordaSec}

The 18th-century French mathematician Jean-Charles de Borda devised
several counting mechanisms. We elaborate one such method, often
called the \emph{modified} De Borda count --- but we drop this
qualification. This counting method can be described abstractly and
efficiently using that multisets form a \emph{monad}. This means that
it comes with two special operations, which we call unit and
flatten. They satisfy certain equations, which are not directly
relevant here.

\begin{definition}
\label{MltMonadDef}
For an arbitrary set $X$ there are functions $\unit \colon X
\rightarrow \natMlt(X)$ and $\flatten \colon
\natMlt\big(\natMlt(X)\big) \rightarrow \natMlt(X)$ given as:
\[ \begin{array}{rclcrcl}
\unit(x)
& = &
1\ket{x}
& \qquad &
\flatten\Big(\sum_{i}\, n_{i}\Bigket{\varphi_{i}}\Big)
& = &
\sum_{i}\, n_{i}\cdot \varphi_{i}.
\end{array} \]

\noindent This last expression uses (repated) addition on multisets,
where $n\cdot \varphi = \varphi + \cdots + \varphi$.
\end{definition}

The unit map turns an element into a singleton multiset. That is easy.
The flatten map `flattens' a multiset of multisets into a multiset.
What actually happens is not so complicated, as illustrated in:
\[ \begin{array}{rcl}
\lefteqn{\flatten\Big(2\Bigket{3\ket{a} + 4\ket{b}} + 3\Bigket{5\ket{a} + 6\ket{b}}\Big)}
\\[+0.4em]
& = &
\Big(3\ket{a} + 4\ket{b}\Big) + \Big(3\ket{a} + 4\ket{b}\Big) +
   \Big(5\ket{a} + 6\ket{b}\Big) + \Big(5\ket{a} + 6\ket{b}\Big) +
   \Big(5\ket{a} + 6\ket{b}\Big)
\\[+0.2em]
& = &
21\ket{a} + 26\ket{b}.
\end{array} \]

We explain the De Borda count also via an example. We actually reuse
the situation described in Table~\eqref{PrefSeqTable}. A preference
sequence $\tuple{x_{1}, x_{2}, x_{3}, x_{4}}$ is now counted in a
weighted manner, like in Example~\ref{ThreeOptionEx}. Let's say the
first preference counts for four, the second for three, the third for
two and the fourth for one. This means that we interprete the
preference sequence $\tuple{x_{1}, x_{2}, x_{3}, x_{4}}$ as a multiset
$4\ket{x_{1}} + 3\ket{x_{2}} + 2\ket{x_{3}} + 1\ket{x_{4}}$. The whole
table in~\eqref{PrefSeqTable} is then interpreted as a multiset over
multisets. The election outcome is now obtained via flattenening.
Explicitly:
\[ \begin{array}{rcl}
\lefteqn{\flatten\Big(7\Bigket{4\ket{a} \shortplus 3\ket{b} \shortplus 2\ket{d} \shortplus 1\ket{c}} + 
   2\Bigket{4\ket{a} \shortplus 3\ket{c} \shortplus 2\ket{d} \shortplus 1\ket{b}}}
\\[+0.2em]
& & \quad \,+
   7\Bigket{4\ket{a} \shortplus 3\ket{d} \shortplus 2\ket{c} \shortplus 1\ket{b}} + 
   8\Bigket{4\ket{b} \shortplus 3\ket{a} \shortplus 2\ket{c} \shortplus 1\ket{d}}
\\[+0.2em]
& & \quad + \,
   7\Bigket{4\ket{b} \shortplus 3\ket{c} \shortplus 2\ket{d} \shortplus 1\ket{a}} + 
   19\Bigket{4\ket{b} \shortplus 3\ket{d} \shortplus 2\ket{c} \shortplus 1\ket{a}} 
\\[+0.2em]
& & \quad + \,
   12\Bigket{4\ket{c} \shortplus 3\ket{a} \shortplus 2\ket{d} \shortplus 1\ket{b}} + 
   12\Bigket{4\ket{d} \shortplus 3\ket{a} \shortplus 2\ket{b} \shortplus 1\ket{c}}
\\[+0.2em]
& & \quad + \,
   8\Bigket{4\ket{d} \shortplus 3\ket{a} \shortplus 2\ket{c} \shortplus 1\ket{b}} + 
   6\Bigket{4\ket{d} \shortplus 3\ket{b} \shortplus 2\ket{c} \shortplus 1\ket{a}}
\\[+0.2em]
& & \quad + \,
   12\Bigket{4\ket{d} \shortplus 3\ket{c} \shortplus 2\ket{a} \shortplus 1\ket{b}}\Big)
\\[+0.2em]
& = &
240\ket{a} + 240\ket{b} + 226\ket{c} + 294\ket{d}.
\end{array} \]

\noindent Candidate $d$ thus wins.

In this example we have mapped a sequence $\tuple{x_{1}, x_{2}, x_{3},
  x_{4}}$ to a multiset $4\ket{x_{1}} + 3\ket{x_{2}} + 2\ket{x_{3}} +
1\ket{x_{4}}$. There are alternative ways to do this. For instance, in
the Eurovision Song contest, each participating country votes by
providing a list of ten (other) countries $\tuple{x_{1}, \ldots,
  x_{10}}$.  They are not counted as $10\ket{x_1} + 9\ket{x_2} +
\cdots + 1\ket{x_{10}}$, but as $12\ket{x_1} + 10\ket{x_2} + 8\ket{x_3} +
7\ket{x_2} + \cdots + 1\ket{x_{10}}$. This gives rise to the famous
phrase: \emph{the twelve points go to \ldots} The next time you hear
it you may wish to think of the monad at work.

\section{Counting single transferable votes}\label{STVSec}

We assume again in this section that a voter can vote 
by numbering candidates in the order of his or her preference 
where, for simplicity, we assume that the resultant ordered 
lists of candidates all have the same length, the length being 
equal to the number of candidates (full preferential voting). 

All the vote counting presented so far has involved finding a
`winner' --- a single successful candidate or result.  Another 
important kind of voting system elects multiple members from the 
one electorate, aiming to better represent the interests of the 
community as a whole rather than just of the majority as measured 
by some vote counting technique.

In such a multi-member electoral situation, the goal is not to find a
single winner, as in Section~\ref{InstantSec}, but to fill multiple
seats --- so we need to determine multiple winners. If there are $N$
votes and $S$ seats, then a candidate with greater than $\frac{N}{S}$
votes deserves a seat. This fraction is called the `quota'. It is
computed below in a slightly different manner, in so-called Droop
style.

A candidate may well receive many more votes than the quota. Single
transferable vote counting (STV) aims to ensure that
these `surplus' votes are not lost, but are redistributed 
(`transferred') according to the next preference on the list
after the chosen candidate. 

It is not a priori obvious how to arrange the surplus transfer.
Suppose that the quota is $Q$ and that $K > Q$ votes are cast in 
favour of candidate $k$.  Various methods have
been tried including ``last parcel'' (transfer all, and only, 
the $K-Q$ votes which occur in the count process after the 
first $Q$ votes for candidate $k$ are counted) and ``random selection'' 
(choose randomly from all $K$ votes for candidate $k$ a selection
of $K-Q$ votes to transfer).  Both of those techniques unfortunately
mean that the outcome can depend upon the order in which votes are
counted or randomly chosen.  A better approach is to transfer all
$K$ votes to the next preferred candidate, but with fractional values
equal to $\frac{K-Q}{K}$, called  the `transfer
value' fraction.  Equivalently, and this is the approach we illustrate
below, the $K$ votes for candidate $k$ are devalued, once $k$ is 
considered elected, by $Q/K$, the fact that they were used to elect
candidate $k$ is thus taken into account, and they are then used with 
that reduced value in favour of the next preferred candidate.

The transfer process itself needs to be
iterated as more and more seats are filled, 
and again various methods have been tried.
Fortunately the multiset
perspective steers us in the right direction.

The process will be worked through below.  We illustrate the so-called
`weighted Gregory method'~\cite{FarrellM06} which is closely related
to the vote counting in the Western Australian upper
house\footnote{For further details, see the website of the
\href{https://www.parliament.wa.gov.au/WebCMS/WebCMS.nsf/resources/file-31-proportional-representation-lc/$file/Sheet\%2031\%20-\%20Proportional\%20Representation\%20in\%20the\%20Legislative\%20Council.pdf}{Western
  Australian Electoral Commission}.}.

Because the transfer value calculation leads to
fractional votes that will be redistributed, 
we need to relax the concept of multiset and allow fractions
instead of only natural numbers as multiplicities --- as in
Notation~\ref{NatMltNot}.

\begin{notation}
\label{MltNot}
For a set $X$, we write $\Mlt(X)$ for the set of multisets
$q_{1}\ket{x_1} + \cdots + q_{K}\ket{x_K}$, with elements $x_{i}\in X$
and with non-negative fractions $q_{i}\in\mathbb{Q}_{\geq 0}$ as
multiplicities. This $\Mlt$ is also a functor, and even a monad. 

%% In fact, we could be more general and allow $q_{i}$ to be non-negative
%% real numbers.
\end{notation}

There is another new element that we use below, namely subtraction
$\psi - \varphi$ of multisets to do the 'devaluation' by $Q/K$ 
described above. This works in the obvious way,
elementwise, exactly like multiset addition. We ensure 
of course that no negative multiplicities arise.

For the example of single transferable votes, we assume that there are
five candidates $a,b,c,d,e$ and a table with 250 vote sequences of
the following form.
\begin{equation}
\label{PrefSeqSTVTable}
\hbox{\begin{tabular}{c||c|c|c|c|c|c|c|c|c|c|c|c|c|c|c|c}
 & 20 & 23 & 27 & 23 & 
   7 & 14 & 13 & 17 & 
   12 & 13 & 7 & 14 & 
   9 & 15 & 16 & 20
\\
\hline
\hline
$1^{\mathsl{st}}$ & 
   $a$ & $a$ & $a$ & $a$ & 
   $b$ & $b$ & $c$ & $c$ & 
   $c$ & $d$ & $d$ & $d$ & 
   $e$ & $e$ & $e$ & $e$
\\
$2^{\mathsl{nd}}$ & 
   $b$ & $b$ & $c$ & $e$ & 
   $a$ & $e$ & $a$ & $d$ & 
   $e$ & $b$ & $c$ & $e$ & 
   $a$ & $b$ & $c$ & $d$
\\
$3^{\mathsl{rd}}$ & 
   $c$ & $d$ & $e$ & $b$ & 
   $d$ & $a$ & $d$ & $b$ & 
   $a$ & $c$ & $b$ & $a$ & 
   $c$ & $c$ & $b$ & $b$
\\
$4^{\mathsl{th}}$ & 
   $e$ & $e$ & $d$ & $c$ & 
   $e$ & $d$ & $e$ & $a$ & 
   $d$ & $e$ & $e$ & $c$ & 
   $b$ & $d$ & $a$ & $a$
\\
$5^{\mathsl{th}}$ & 
   $d$ & $c$ & $b$ & $d$ & 
   $c$ & $c$ & $b$ & $e$ & 
   $b$ & $a$ & $a$ & $b$ & 
   $d$ & $a$ & $d$ & $c$
\end{tabular}}
\end{equation}

\noindent This table translates into the following multiset
$\varphi_{0} \in \natMlt\big(\Lst(\{a,b,c,d,e\})\big) \subseteq
\Mlt\big(\Lst(\{a,b,c,d,e\})\big)$.
\[ \begin{array}{rcl}
\varphi_{0}
& = &
20\bigket{a, b, c, e, d} + 
23\bigket{a, b, d, e, c} + 
27\bigket{a, c, e, d, b} + 
23\bigket{a, e, b, c, d} 
\\[+0.2em]
& & \quad + \, 
7\bigket{b, a, d, e, c} +
14\bigket{b, e, a, d, c} + 
13\bigket{c, a, d, e, b} + 
17\bigket{c, d, b, a, e} 
\\[+0.2em]
& & \quad + \, 
12\bigket{c, e, a, d, b} +
13\bigket{d, b, c, e, a} +
7\bigket{d, c, b, e, a} + 
14\bigket{d, e, a, c, b} 
\\[+0.2em]
& & \quad + \, 
9\bigket{e, a, c, b, d} + 
15\bigket{e, b, c, d, a} + 
16\bigket{e, c, b, a, d} + 
20\bigket{e, d, b, a, c}.
\end{array} \]

\noindent In this illustration we assume that there are $S_{0}=3$
seats to be filled. The Droop quota function is defined, on arbitrary
arguments, as:
\[ \begin{array}{rcl}
\quota(\varphi, S)
& \coloneqq &
\displaystyle \left\lfloor \frac{\|\varphi\|}{S+1} \right\rfloor + 1.
\end{array} \]

The Droop quota, which is generally smaller than 
$\frac{N}{S} = \frac{\|\varphi\|}{S}$,
is used because if $S$ candidates achieve greater than or equal to 
the Droop quota, then no other candidate can possibly receive more 
votes than each of them, and so those $S$ candidates will be 
declared elected.

We go through the following iterative steps.
\begin{enumerate}
\item Starting form the multiset $\varphi_{0}$ and the number of seats
  $S_{0}=3$ we compute $\quota\big(\varphi_{0}, S_{0}\big) = \lfloor
  \frac{250}{4} \rfloor + 1 = \lfloor 62.5 \rfloor + 1 = 63$. We then
  check if there is a candidate who has at least 63 votes as first
  preference. We do this, like in Section~\ref{InstantSec}, via
  functoriality of $\Mlt$ applied to the the head-of-list function
  $\head$.  This gives:
\[ \begin{array}{rcl}
\Mlt\big(\head\big)(\varphi_{0})
& \smash{\stackrel{\eqref{MltFunEqn}}{=}} &
93\ket{a} + 21\ket{b} + 42\ket{c} + 34\ket{d} + 60\ket{e}.
% 93|a> + 21|b> + 42|c> + 34|d> + 60|e>
\end{array} \]

\noindent We see that only candidate $a$ has more than $63$ first
preference votes, and thus gets a seat. There are
$93-63 = 30$ surplus votes for $a$ that need to be
redistributed. Formulated differently, we need to subtract from
$\varphi_{0}$ the fraction $\frac{Q}{K} = \frac{63}{93}$ 
of the sequences that start
with $a$. This fraction of the votes was used up for the seat
of~$a$. This subtraction takes the form:
\[ \varphi_{0} - \frac{63}{93} \cdot 
   \Big(20\bigket{a, b, c, e, d} + 23\bigket{a, b, d, e, c} + 
   27\bigket{a, c, e, d, b} + 23\bigket{a, e, b, c, d}\Big). \]

\noindent We obtain the new multiset of votes $\varphi_{1}$ by deleting
candidate $a$ from this difference:
\[ \begin{array}{rcl}
\varphi_{1}
& \coloneqq &
\displaystyle\Mlt\big(\delete_{a}\big)\left(
   \varphi_{0} - \frac{63}{93} \cdot 
   \Big(20\bigket{a, b, c, e, d} + 23\bigket{a, b, d, e, c} + 
   27\bigket{a, c, e, d, b} + 23\bigket{a, e, b, c, d}\Big)\right)
\\[+0.4em]
& = &
6.45\bigket{b, c, e, d} + 
14.42\bigket{b, d, e, c} + 
14\bigket{b, e, d, c} + 
17\bigket{c, d, b, e} + 
\\[+0.2em]
& & \quad + \, 
13\bigket{c, d, e, b} + 
20.71\bigket{c, e, d, b} + 
13\bigket{d, b, c, e} + 
7\bigket{d, c, b, e} + 
\\[+0.2em]
& & \quad + \, 
14\bigket{d, e, c, b} + 
22.42\bigket{e, b, c, d} + 
25\bigket{e, c, b, d} + 
20\bigket{e, d, b, c}.
\end{array} \]

% decimals: 42 + 71 + 42 + 45 = 200

\noindent One sees that in this way candidate $b$ gets quite a few
extra votes because it is listed as second in $20 + 23 = 43$ cases
where people had $a$ as first preference. Of all of those, a
`transfer' fraction of $1 - \frac{63}{93} = \frac{30}{93}$ with $b$ as
first preference appears in $\varphi_{1}$. Candidates $c$ and $e$ also
get additional (first preference) votes, but not $d$.

The new multiset $\varphi_{1}$ contains fractions as multiplicities,
and thus forms an element of the set
$\Mlt\big(\Lst(\{a,b,c,d,e\})\big)$, see Notation~\ref{MltNot}. One
may check that the size of the multiset $\varphi_{1}$ is a natural
number, namely $250-63 = 187$. This is the number of (as yet)
unused votes. The number of available seats at this stage is $S_{1}
\coloneqq S_{0} - 1 = 2$.

In our example, only one candidate has reached the quota. There could
be more. In that case, for each winning candidate $x$, the sequences
starting with $x$ are subtracted from $\varphi_{0}$, with
corresponding fractions given by the quota divided by the number of
first preference votes for $x$.

\item In the next step we (happen to) have the same quota, but now
  computed as $\quota\big(\varphi_{1}, S_{1}\big) = \lfloor
  \frac{187}{3} \rfloor + 1 = \lfloor 62.333 \rfloor + 1 = 63$. The
  first preferential votes are at this stage:
\[ \begin{array}{rcl}
\Mlt\big(\head\big)(\varphi_{1})
& = &
34.87\ket{b} + 50.71\ket{c} + 34\ket{d} + 67.42\ket{e}.
% 34.87|b> + 50.71|c> + 34|d> + 67.42|e>
\end{array} \]

\noindent We see that candidate $e$ now wins a seat, with a small
surplus of $67.42 - 63$. We adapt the number of available seats to
$S_{2} \coloneqq S_{1} - 1 = 1$, and we form the next successor
multiset of votes:
\[ \begin{array}{rcl}
\varphi_{2}
& \coloneqq &
\displaystyle\Mlt\big(\delete_{e}\big)\left(
   \varphi_{1} - \frac{63}{67.42} \cdot 
   \Big(22.42\bigket{e, b, c, d} + 25\bigket{e, c, b, d} + 
   20\bigket{e, d, b, c}\Big)\right)
\\[+0.8em]
& = &
7.92\bigket{b, c, d} + 
28.42\bigket{b, d, c} + 
1.64\bigket{c, b, d} + 
50.71\bigket{c, d, b} + 
14.31\bigket{d, b, c} + 
21\bigket{d, c, b}.
\end{array} \]

%% 7.921|['b', 'c', 'd']> + 
%% 28.42|['b', 'd', 'c']> + 
%% 1.639|['c', 'b', 'd']> + 
%% 50.71|['c', 'd', 'b']> + 
%% 14.31|['d', 'b', 'c']> + 
%% 21|['d', 'c', 'b']>

% decimals 92 + 42 + 64 + 71 + 31 = 300

\noindent This multiset has size $250 - 63 - 63 = 124$.

\item In the next round we get again as quota:
  $\quota\big(\varphi_{2}, S_{2}\big) = \lfloor \frac{124}{2} \rfloor
  + 1 = \lfloor 62 \rfloor + 1 = 63$. The first preferences are at
  this stage:
\[ \begin{array}{rcl}
\Mlt\big(\head\big)(\varphi_{2})
& = &
36.34\ket{b} + 52.35\ket{c} + 35.31\ket{d}.
%36.34|b> + 52.35|c> + 35.31|d>
\end{array} \]

\noindent There is no winner. Then, candidate $b$ with the lowest number
of votes is removed from the current multiset $\varphi_{2}$, giving 
as next multiset:
\[ \begin{array}{rcccl}
\varphi_{3}
& \coloneqq &
\displaystyle\Mlt\big(\delete_{b}\big)\big(\varphi_{2}\big)
& = &
60.27\bigket{c, d} + 63.73\bigket{d, c}.
% 60.27|['c', 'd']> + 63.73|['d', 'c']>
\end{array} \]

\item We now trivially get:
\[ \begin{array}{rcl}
\Mlt\big(\head\big)(\varphi_{3})
& = &
60.27\bigket{c} + 63.73\bigket{d},
\end{array} \]

\noindent with quota $\quota\big(\varphi_{3}, S_{2}\big) = 63$, so
that candidate $d$ gets the last seat.
\end{enumerate}

\section{Conclusions}\label{ConclusionSec}

This paper illustrates (and demonstrates) the usefulness of multisets,
and of their categorical properties, in vote counting algorithms that
are used in practice. This is applied category theory in elections.
The paper is not `deep' from a categorical perspective, but it
introduces categorical techniques in a new field where they have not
been identified before.

We believe that this usage of multisets brings clarity to an area
where transparency and correctness are important. There is more work
to do, both at a theoretical level and at a more practical level.
\begin{itemize}
\item This paper concentrates on representing the counting algorithms
  in terms of multisets. This multiset representation can now be used
  to prove basic properties about these algorithms. For instance, a
  first requirement is stability of the outcome under permutation of
  votes. This is automatic when one uses multisets, since the order of
  their elements (votes) is irrelevant.

\item We have used various simplifying assumptions that may not hold
  in practice, such as: all voters express lists of candidate
  preferences of the same length. Our descriptions can be adapted to
  situations where this fails, but they then have to deal with
  `undefinedness', for instance when using the head-of-list
  operation. This is manageable, but makes the descriptions less
  smooth. Edge cases like a tied vote are also not covered here.

\item Given the new abstract vote count algorithms presented here, it
  makes sense to develop reference implementations of the frequently
  used versions in programming languages that support multisets, like
  Haskell. One could then do a recount or a parallel count of earlier
  and upcoming elections. And of course, these multiset-based
  implementations could be used in actual, future elections.

% https://hackage.haskell.org/package/multiset-0.3.4.3/docs/Data-MultiSet.html

\end{itemize}

Historically, vote counting was based on sorting and systematically
transferring physical ballot papers.  Most modern elections now 
use electronic voting or scanned or manually entered ballot papers
to permit computer supported counting.  The multiset approach to 
vote counting liberates the algorithms from their old association 
with physical ballot papers, focusing instead on common sequences
of preferences (the lists that are the basic elements of the 
multisets).  The multiset approach, and the categorical properties 
of multisets, significantly simplify the vote counting process
with potential benefits for security, accuracy and explainability.

%% Further comments of Mike

%% By the way, you know I guess that Gregory was Australian?  (Or at
%% least Australians say so -- I haven't checked carefully and we do
%% sometimes claim New Zealanders or Canadians if they live here long
%% enough!)

%\bibliography{multiset-counting.bib}

%\bibliography{/home/bart/svn/bart/Tex/bib}

\end{document}